\documentclass[journal,transmag]{IEEEtran}
\ifCLASSINFOpdf
  \usepackage[pdftex]{graphicx}
\else
\fi
\usepackage{amsmath}
\usepackage{amssymb}
\newcommand{\opcurl}{\operatorname{curl}}
\newcommand{\MR}{\mathbb R}

\begin{document}
%
\title{An Iterative Method for Transient Finite Element Simulations of Non-Linear Eddy Current Problems}


\author{\IEEEauthorblockN{K. Hollaus\IEEEauthorrefmark{1},
and H. Silm\IEEEauthorrefmark{1}}
\IEEEauthorblockA{\IEEEauthorrefmark{1}Institute of Analysis and Scientific Computing, TU Wien, Vienna, A-1040 Austria}
}

%



\IEEEtitleabstractindextext{%
\begin{abstract}
A method is presented to carry out a transient simulation of eddy current problems with nonlinear materials. Coils are voltage-driven. The magnetic field due to currents in coils are considered by their Biot-Savart-fields. The magnetic vector potential is used in the finite element formulation. The time stepping method is based on implicit Euler. The arising nonlinear equation system is split into two parts, the common finite element system and a circuit equation. Each part is solved separately by Newton's method. Additionally, a line search is used to solve the nonlinear field equations. Inrush currents and average magnetic flux densities through cross sections of laminates of a nonlinear benchmark problem consisting of a laminated iron core inserted in a cylindrical coil are studied. All details of the numerical benchmark are given to evaluate the presented results. Numerical data describing the performance of the presented method are provided. 
\end{abstract}

}

\maketitle

\IEEEdisplaynontitleabstractindextext

%
\IEEEpeerreviewmaketitle

\section{Introduction}
%
%
%
%
\IEEEPARstart{T}{he method} presented here facilitates a transient finite element simulation of nonlinear eddy current problems with the magnetic vector potential (MVP) $\boldsymbol{A}$ taking into account of voltage-driven coils, of currents in coils by Biot-Savart-field (BSF) and of material properties by a magnetization curve. Such a method is required for instance to study inrush currents of electrical devices. Inrush currents are of great importance in the design of electrical devices due to the electrical and mechanical stability. Reference solutions may also support the development of approximate methods considering laminated iron cores efficiently.

An early work of a coupling of a voltage-driven coil with the finite element method (FEM) based on $\boldsymbol{A}$ can be found in \cite{LeoRod:88}. A method based on a current vector potential $\boldsymbol{T}$ and a reduced magnetic scalar potential $\Phi$ is presented in \cite{MeunFloGuer03}. It can handle solid conductors. A comprehensive overview of various possible coupling techniques for different problems can be found in the work \cite{DULAR:04}. The solution of the arising non-linear system is not discussed by the above works. \\
A coupling of the finite element harmonic balance method (FEHBM) using $\boldsymbol{T}$ and $\Phi$ with a circuit to get the steady state solution has been introduced in \cite{PlaKoczBir:16}. The nonlinear system is iteratively solved with the aid of a fixed point method technique allowing a parallel processing of the resulting problem. The problem is solved combined and separated.

Iterative solvers are required for large linear equation systems. To get a reasonable convergence rate pre-conditioners are required. The problem arises, what is the appropriate pre-conditioner for the present case? The system is composed of two parts, the common non-linear system due to the FEM with $\boldsymbol{A}$ extended by the non-linear equation due to the circuit coupling. To circumvent this problem the system has been split corresponding to the two parts. 

The arising nonlinear equation system is split into two parts, Newton's method (NM) is applied to each part separately exploiting the fast convergence rate of NM. The two parts are alternately solved until a stopping criterion is fulfilled. NM for the system due to the FEM is supplemented by a line search. The BSF caused by a current in a coil is calculated only once and then accordingly scaled to the course of time. 
Induction effects in the windings are neglected (stranded coils). A solution for that can be found for linear materials and BSF in \cite{BiPrBuTi:04}. 

The paper presents the case of one coil in detail, the extension to several coils is straight forward and outlined at the end. A study of a numerical benchmark is presented. To evaluate the results the geometry with all dimensions, the magnetization curve and the excitation are given.
\setlength{\tabcolsep}{4pt} 
\section{Eddy Current Problem}
\subsection{Boundary Value Problem} 
\begin{figure}[h]
	\begin{center}
		\includegraphics[height=4.5cm]{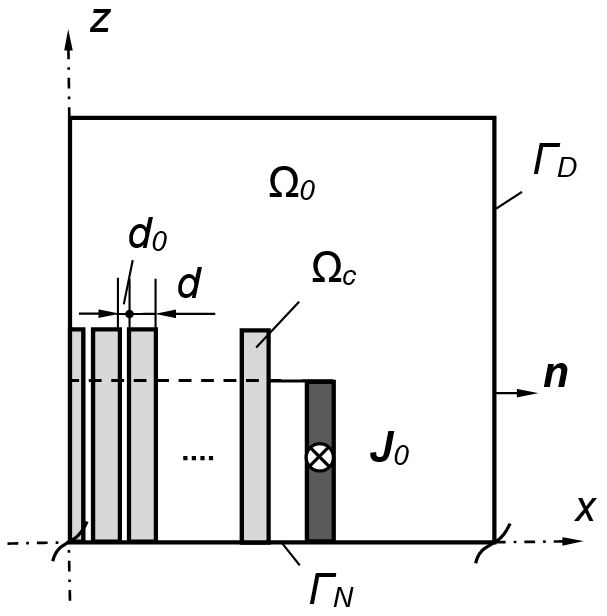}
		\caption{One eighth of the boundary value problem, coil with laminated iron core, $d_1$ and $d_0$ are the thickness of the laminates and the width of the air gap in between.}
		\label{fig:BVP}
	\end{center}
\end{figure}

The eddy current problem to be solved is sketched in Fig.~\ref{fig:BVP}. It consists of a conducting domain (iron) $\Omega_c$ and air $\Omega_0$, i.e., $\Omega = \Omega_c \cup \Omega_0$ with the boundary $\Gamma = \Gamma_D \cup \Gamma_N$. The boundary value problem with $\boldsymbol{A}$ in the time domain reads as  
\begin{eqnarray}
		\opcurl \frac{1}{\mu(\boldsymbol{A})} \opcurl \boldsymbol{A} + \sigma \frac{\partial}{\partial t} \boldsymbol{A} & \ = 
			\boldsymbol{J}_0 & \text{in} \ \Omega \subset \MR ^3, \label{eq1} \\
		\boldsymbol{A} \times \boldsymbol{n} & = \boldsymbol{\alpha} & \text{on} \ \Gamma_D, \label{eq2} \\
		\frac{1}{\mu(\boldsymbol{A})} \opcurl \boldsymbol{A} \times \boldsymbol{n} & \ = \boldsymbol{K} & \text{on} \ \Gamma_N, \label{eq3} 
\end{eqnarray}
where $\boldsymbol{J}_0$ in (\ref{eq1}) stands for a source current density in $\Omega_s \subset \Omega_0$ known at the moment, $\boldsymbol{\alpha}$ in (\ref{eq2}) represents a magnetic flux and $\boldsymbol{K}$ in (\ref{eq3}) describes a surface current density. Both, $\boldsymbol{\alpha}$ and $\boldsymbol{K}$, are zero for the considered problem. Indices $D$ and $N$ stand for Dirichlet and Neumann boundary, respectively. The material parameters are the magnetic permeability $\mu(\boldsymbol{A})$, which is assumed to be nonlinear, and the electric conductivity $\sigma$, respectively. \\

\subsection{Weak Form} 
Equations (\ref{eq1}) to (\ref{eq3}) lead to the following weak form for the FEM.
Find $\boldsymbol{A}_h \in \mathcal{V}_{\alpha} := \{\boldsymbol{A}_h \in \mathcal{V}_h: \boldsymbol{A}_h  \times  \boldsymbol{n} = 
{\boldsymbol{\alpha}}_h$ on $\Gamma \}$, such that
\begin{eqnarray}
	\int_{\Omega} \frac{1}{\mu(\boldsymbol{A}_h)} \opcurl \boldsymbol{A}_h \cdot \opcurl \boldsymbol{v}_h \,d\Omega +  \label{eq4}  \\
	\frac{\partial}{\partial t} \int_{\Omega}\sigma \boldsymbol{A}_h \cdot  \boldsymbol{v}_h \,d\Omega =
	\int_{\Omega_s} \boldsymbol{J}_0 \cdot \boldsymbol{v}_h \,d\Omega \nonumber 
\end{eqnarray} 
for all $\boldsymbol{v}_h \in \mathcal{V}_0 := \{\boldsymbol{v}_h \in \mathcal{V}_h: {\boldsymbol{v}}_h  \times  \boldsymbol{n} = 
\boldsymbol{0}$ on $\Gamma_D \}$, where $\mathcal{V}_h$ is a finite element subspace of $H(\opcurl,\Omega)$. \\
Regularization is applied in $\Omega_0$ to obtain a unique solution \cite{LedgZagl:10}.
\subsection{Voltage-Driven Coil and Biot-Savart-Field} 
In case of a voltage-driven coil and non-linear materials the current density $\boldsymbol{J}_0$ can no longer be presumed. The network equation
\begin{equation}
	u_0(t) = i(t)R_s - \frac{d \phi(t)}{d t} \label{eq5}
\end{equation}
with voltage $u_0$ of the voltage source, the net current $i$, the series resistor $R_s$ representing the resistance of the coil plus a possible series resistor and the magnetic flux $\phi$ has to be solved together with (\ref{eq4}). 

Rearranging of the linear form in (\ref{eq4}) yields
\begin{eqnarray}
	\int_{\Omega_s} \boldsymbol{J}_0 \cdot \boldsymbol{v}_h \,d\Omega = 
	i(t) \int_{\Omega_s} \boldsymbol{\tau} \cdot \boldsymbol{v}_h \,d\Omega =  \nonumber \\
	\int_{\Omega_s} \opcurl \boldsymbol{H}_0 \cdot \boldsymbol{v}_h \,d\Omega = 
	i(t) \int_{\Omega_s} \opcurl \boldsymbol{h}_0 \cdot \boldsymbol{v}_h \,d\Omega,	\label{eq6}
\end{eqnarray}
where $\boldsymbol{\tau}$ is the turn density, $\boldsymbol{H}_0$ is the BSF of $\boldsymbol{J}_0$, $\boldsymbol{j}_0$ is the current density due to a unit net current \cite{DULAR:04}. The BSF $\boldsymbol{h}_0$ of this unit net current is obtained by 
\begin{equation}
	\boldsymbol{h}_{0} = \frac{1}{4 \pi} \int_{\Omega_s} \frac{\boldsymbol{j}_0 \times (\boldsymbol{r}_F-\boldsymbol{r}_S)}
		{ | \boldsymbol{r}_F-\boldsymbol{r}_S | ^3} \,d\Omega,   \label{eq7}
\end{equation}
where $\boldsymbol{r}_F$ and $\boldsymbol{r}_S$ are the field and the source point.
The integral in the last term of (\ref{eq6}) extended by zero in $\Omega \setminus \Omega_0$ and integrated by parts leads to
\begin{eqnarray}
\int_{\Omega} \opcurl \boldsymbol{h}_0 \cdot \boldsymbol{v}_h \,d\Omega = \int_{\Omega} \boldsymbol{h}_0 \cdot \opcurl \boldsymbol{v}_h \,d\Omega + \nonumber \\
\oint_{\Gamma} (\boldsymbol{h}_0 \times \boldsymbol{v}_h) \cdot \boldsymbol{n} \,d\Gamma.  \label{eq8}
\end{eqnarray}
The surface integral in (\ref{eq8}) vanishes because either $\boldsymbol{h}_0$ or $\boldsymbol{v}_h$ equals to zero on planes of symmetry or $\boldsymbol{h}_0$ is assumed to be zero on the far boundary \cite{BiPrBuTi:04}. Note, that thanks to $\boldsymbol{\tau} = \opcurl \boldsymbol{h}_0$ and that $\boldsymbol{h}_0$ is sufficiently smooth in $\Omega$, respectively, the volume of the coil need not to be modeled by FEs and (\ref{eq8}) can be integrated numerically without any difficulties, respectively.

At the same time rearranging of the second term on the right hand side of (\ref{eq5}) leads to
\begin{eqnarray}
	- \frac{d \phi(t)}{d t} = \int_{S w} \boldsymbol{E} \cdot \,d\boldsymbol{s} = \int_{S} \boldsymbol{E} \cdot \frac{w \boldsymbol{e}_F}{F} \cdot F \boldsymbol{e}_F  \cdot \,d	\boldsymbol{s} =  \nonumber  \\
	\int_{S} \boldsymbol{E} \cdot \boldsymbol{\tau} \cdot F \boldsymbol{e}_F  \cdot \,d\boldsymbol{s} =
	\int_{\Omega_s} \boldsymbol{E} \cdot \boldsymbol{\tau} \,d\Omega = \nonumber \\
	\int_{\Omega_s} \boldsymbol{E} \cdot \opcurl \boldsymbol{h}_0 \,d\Omega =
	\int_{\Omega} \boldsymbol{E} \cdot \opcurl \boldsymbol{h}_0 \,d\Omega = \nonumber \\ 
	\int_{\Omega} \opcurl \boldsymbol{E} \cdot \boldsymbol{h}_0 \,d\Omega - 
	\oint_{\Gamma} (\boldsymbol{E} \times \boldsymbol{h}_0) \cdot \boldsymbol{n} \,d\Gamma = \nonumber  \\
 -\int_{\Omega} \opcurl \frac{\partial \boldsymbol{A}}{\partial t} \cdot \boldsymbol{h}_0 \,d\Omega \label{eq9}
\end{eqnarray}
with $F$, $\boldsymbol{E}$, $S$, $w$ and $\boldsymbol{e}_F$, respectively, the cross-section of the coil, the electric field strength, the average path length of a winding, the number of turns and the unit vector of the surface $F$, respectively. Similar considerations have been made in (\ref{eq9}) for the extension by zero and for the surface integral as in (\ref{eq8}).
\subsection{Time Stepping Scheme} 
Evaluation of (\ref{eq4}) and considering (\ref{eq8}) on the one hand and (\ref{eq5}) with (\ref{eq9}) on the other results in the nonlinear ordinary differential equation system
\begin{eqnarray}
	\boldsymbol{S}(\boldsymbol{a})\boldsymbol{a}(t) + \boldsymbol{M} \frac{d}{d t} \boldsymbol{a}(t) - \boldsymbol{b} i(t) & = & \boldsymbol{0} \label{eq10} \\
	- \boldsymbol{b}^T \frac{d}{d t} \boldsymbol{a}(t) + i(t) R_s & = & u_0(t) \label{eq11}
\end{eqnarray}
In (\ref{eq10}) and (\ref{eq11}), $\boldsymbol{a}$,  $\boldsymbol{b}$, $\boldsymbol{S}$ and $\boldsymbol{M}$, respectively, are the unknown vector of the MVP, the vector according to $\boldsymbol{h}_0$ and (\ref{eq8}), the stiffness matrix and the mass matrix, respectively. 
Applying implicit Euler as time stepping scheme leads to the system of nonlinear equations
\begin{eqnarray}
	\boldsymbol{S}(\boldsymbol{a}^{k+1}) \boldsymbol{a}^{k+1} + \frac{1}{\Delta t} \boldsymbol{M} \boldsymbol{a}^{k+1}- \boldsymbol{b} i^{k+1}  & = & \frac{1}{\Delta t} \boldsymbol{M} \boldsymbol{a}^{k}  \label{eq12} \\
	- \boldsymbol{b}^T \boldsymbol{a}^{k+1} + \Delta t i^{k+1} R_s & = & \Delta t u_0^{k+1} - \boldsymbol{b}^T \boldsymbol{a}^{k} \nonumber \\ 		\label{eq13}
\end{eqnarray}
to be solved, where $\Delta t$ is the time step and $k$ the index for the time instant $t_k = k \Delta t$. The system of (\ref{eq12}) and (\ref{eq13}) is symmetric.

The extension to an arbitrary number $n$ of coils is straightforward. To this end (\ref{eq5}) is replaced by the corresponding  $n$ network equations, (\ref{eq8}) and (\ref{eq9}) are written for each coil resulting in $n$ different vectors $\boldsymbol{b}_i$, $i=1,...,n$.

\subsection{Newton's Method}
The coefficient matrix of the extended equation system (\ref{eq12}) and (\ref{eq13}) does not have the typical structure of a FE equation system anymore. Therefore the existing precondition techniques would have to be modified. To avoid this a separated NM is used. The system to be solved is split into two parts according to (\ref{eq12}) and (\ref{eq13}), i.e. the unknown vector $\boldsymbol{a}$ and the net current $i$.

Let $\boldsymbol{F}(\boldsymbol{a}^{k+1}) = \boldsymbol{0}$ and $G(i^{k+1}) = 0$ be the implicit representations of (\ref{eq12}) and (\ref{eq13}).
Applying NM to (\ref{eq12}) and (\ref{eq13}) yields: 
\begin{eqnarray}
	\boldsymbol{a}^{k+1}_{l+1} & = & \boldsymbol{a}^{k+1}_{l} - \alpha \boldsymbol{J}_F^{-1}(\boldsymbol{a}^{k+1}_{l}) \boldsymbol{F}(\boldsymbol{a}^{k+1}_{l}) \label{eq14} \\
	i^{k+1}_{l+1} & = & i^{k+1}_{l} - J_G^{-1}(i^{k+1}_{l}) G(i^{k+1}_{l}) \label{eq15}
\end{eqnarray}
Index $l$ denotes the nonlinear iterations. The Jacobian matrices are:
\begin{eqnarray}
	\boldsymbol{J}_F(\boldsymbol{a}^{k+1}) & = & \frac{d}{d \boldsymbol{a}^{k+1}} \boldsymbol{F}(\boldsymbol{a}^{k+1})  \label{eq16} \\
	J_G(i^{k+1}) & = & \frac{d}{d i^{k+1}} G(i^{k+1}) \label{eq17}
\end{eqnarray}
To get (\ref{eq17}) 
\begin{equation}
	\boldsymbol{J}_F(\boldsymbol{a}^{k+1}) \frac{d \boldsymbol{a}^{k+1}}{d i^{k+1}} = \boldsymbol{b} \label{eq18}
\end{equation}
has to be solved first. \\
While $\boldsymbol{a}^{k+1}$ of (\ref{eq14}) is solved iteratively the net current $i^{k+1}$ is constant and known and vice versa for (\ref{eq15}).
Stopping criteria are presented in section III. C.
NM to solve (\ref{eq14}) is supplemented by a line search using the functional
\begin{eqnarray}
	J(\boldsymbol{A}) =  \int\limits_{\Omega}\! w(|\opcurl{\boldsymbol{A}}|) \, \mathrm{d}\Omega + \int\limits_{\Omega} \!\frac{\sigma}{2\Delta t}
	\boldsymbol{A}^{\mathrm{T}}\boldsymbol{A}\, \mathrm{d}\Omega - \nonumber \\
	\int\limits_{\Omega}\! \vec{J_0}^{\mathrm{T}} \boldsymbol{A} \, \mathrm{d} 
	\Omega - \int\limits_{\Omega}\! \frac{\sigma}{\Delta t} \boldsymbol{A}_{\mathrm{old}}^{\mathrm{T}}\boldsymbol{A} \, \mathrm{d}\Omega \label{eq19}
\end{eqnarray}
corresponding to (\ref{eq4}) of the eddy current problem to determine the damping factor $\alpha$ in the line search algorithm for the damped NM, to ensure monotonous convergence of the iterative algorithm. The first integral in \eqref{eq19} is the energy stored in the magnetic field with the energy density
\begin{equation}
  w(B) = \int_0^B H(B') \, \mathrm{d}B', \label{eq20}
\end{equation}
where $B=|\opcurl \boldsymbol{A}|$.
\section{Numerical Benchmark}
The correctness of the method has been evaluated by the example presented in \cite{LeoRod:88}, which provides numerical and experimental results.
\subsection{Benchmark Data}
\begin{figure}[h]
	\begin{center}
		\includegraphics[height=3.5cm]{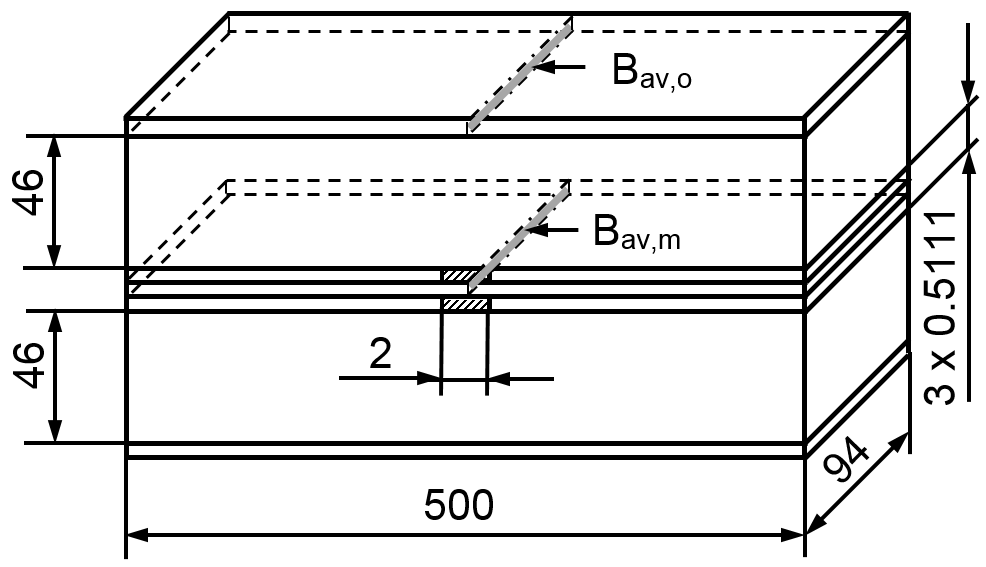}
		\caption{Laminated iron core with dimensions in mm.}
		\label{fig:LIC}
	\end{center}
\end{figure}

\begin{figure}[h]
	\begin{center}
		\includegraphics[height=4.0cm]{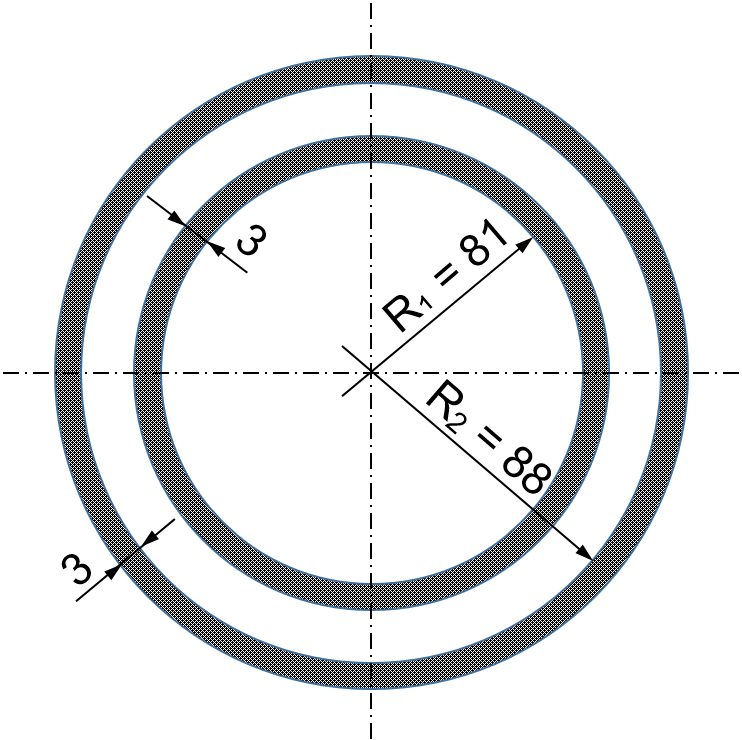}
		\caption{Cross section of the cylindrical coil with dimensions in mm.}
		\label{fig:CCC}
	\end{center}
\end{figure}

The three dimensional numerical benchmark consists of the laminated iron core, see Fig.~\ref{fig:LIC}, inserted in the cylindrical coil shown in Fig.~\ref{fig:CCC} as indicated in Fig.~\ref{fig:BVP}. 
 
The laminated iron core consists of $185$ layers in total and is constructed as follows, see Fig.~\ref{fig:LIC}. From the bottom to the top, the first laminate is followed by 90 laminates. The next layer is composed of $2$ laminates separated by an air gap (hatched, width is $2\,\mathrm{mm}$) in between, next follows one laminate (middle laminate to determine the average flux density $B_{av,m}$ in the cross section $(94\,\mathrm{mm} \times 0.5\,\mathrm{mm})$, grey surface in Fig.~\ref{fig:LIC}), then two separated laminates with an air gap (hatched), then $90$ laminates, and finally $1$ laminate which is used to determine the average flux density $B_{av,o}$ in the cross section (grey surface). 

\begin{figure}[h]
	\begin{center}
		\includegraphics[height=5.0cm]{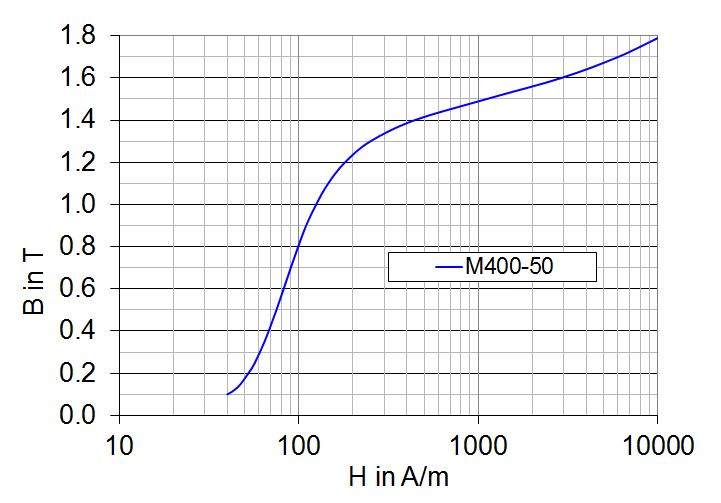}
		\caption{Magnetization curve.}
		\label{fig:MKL}
	\end{center}
\end{figure}

The non-linearity of the isotropic material is described by the magnetization curve in Fig.~\ref{fig:MKL}.

The cross section of the cylindrical coil is shown in Fig.~\ref{fig:CCC}. It consists of two layers (dark rings), $60$ turns per layer. The length of the coil equals $192 \mathrm{mm}$. The frequency was selected with $f = 50 \mathrm{Hz}$. The dimensions are based on a fill-factor of $f_f = 0.9783$, which means an air gap of $d_0=0.01085\,\mathrm{mm}$ between the layers. The coil along with the iron core are arranged symmetrically, i.e. $x=0$, $y=0$ and $z=0$ are planes of symmetry. Therefore, only one eighth of the problem has been modeled by FEs.
Hexahedral edge elements \cite{SchoeZagl:05} were used to simplify the modeling of the laminates by a hand made mesh. The Biot-Savart field was exploited to avoid the modeling of the cylindrical coil. 

\subsection{Simulation Results}

\begin{figure}[h]
	\begin{center}
		\includegraphics[height=5.0cm]{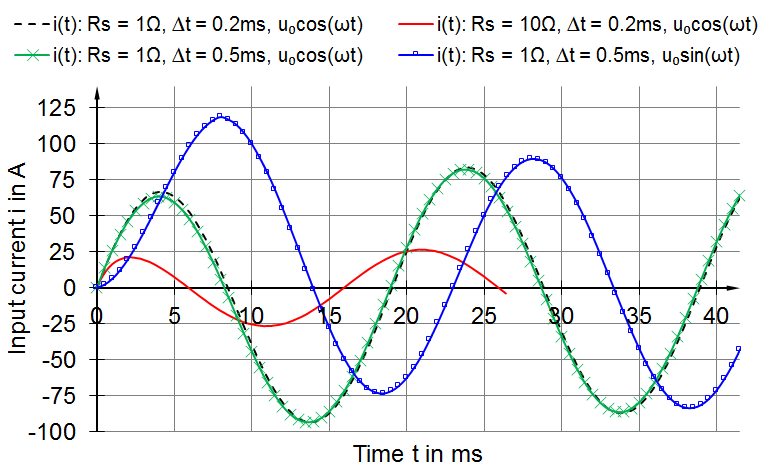}
		\caption{Study of inrush currents for different $R_s$, $\Delta t$ and $u_0(t)$.}
		\label{fig:INCURR}
	\end{center}
\end{figure}

The selected input voltage $u(t)$ equals either $u_0 \cos(\omega t)$ or $u_0 \sin(\omega t)$, where $\omega=2 \pi f$ is the angular frequency and $f$ the frequency. The peak value $u_0$ was assumed to be $280\,\mathrm{V}$ and frequency $50\,\mathrm{Hz}$. Simulations are carried out with $R_s = 1\,\Omega$ and $R_s=10\,\Omega$.
Inrush currents for different parameters are shown in Fig.~\ref{fig:INCURR}.

\begin{figure}[h]
	\begin{center}
		\includegraphics[height=5cm]{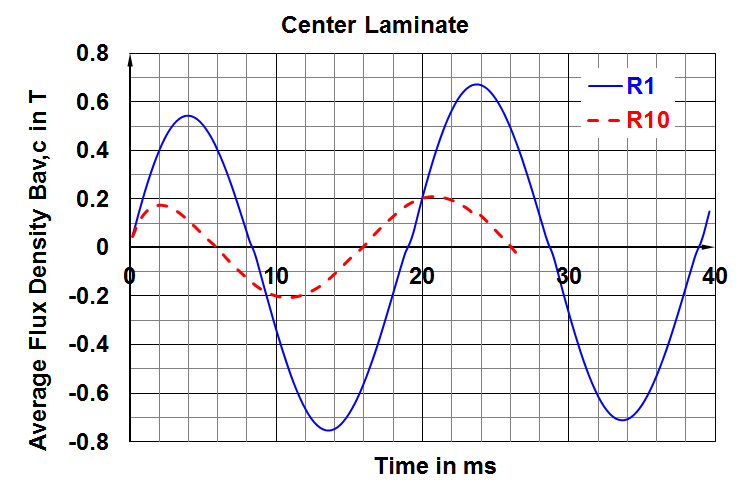}
		\caption{Average flux density in the middle laminate.}
		\label{fig:BAVMI}
	\end{center}
\end{figure}

\begin{figure}[h]
	\begin{center}
		\includegraphics[height=5cm]{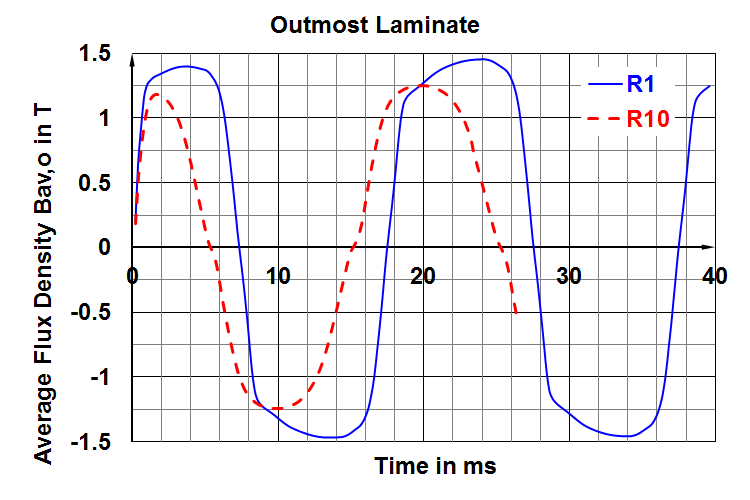}
		\caption{Average flux density in the outmost laminate.}
		\label{fig:BAVOM}
	\end{center}
\end{figure}

The average magnetic flux densities $B_{av}$ for the middle and the out most layer/laminate are shown in Fig.~\ref{fig:BAVMI} and Fig.~\ref{fig:BAVOM}.
The selected time step $\Delta t$ was $0.2\,ms$ and the voltage $u(t)$ was $u_0 \cos(\omega t)$. The solid blue curve belongs to $R = 1.0\,\Omega$ and is marked with $R1$, the dashed red curve to $R = 10.0\,\Omega$ and is marked with $R10$.

Although the currents in Fig.~\ref{fig:INCURR} do not reflect the non-linearity, the magnetic flux density in the out most laminate is clearly nonlinear.
\subsection{Numerical Data}
\begin{center}
Table \, 1: Statistics.\\
\vspace{0.3cm}
\resizebox{\linewidth}{!}{\begin{tabular}{r|ccc|ccc}
  &  & $R=1.0\,\mathrm{\Omega}$ $^a)$ & & & $R=10.0\,\mathrm{\Omega}$ $^a)$ &  \\ \hline 
  & Min. & Max. & Av.  & Min. & Max. & Av. \\ \hline 
  Newton teps FE-system & 2 & 2 & 2 & 2 & 2 & 2 \\ 
  Newton steps net current & 1 & 8 & 3.98 & 1 & 17 & 8.292 \\ 
  Line search steps & 1 & 2 & 1.06 & 1 & 2 & 1.027 
\end{tabular}}
\end{center}
$^a)$ For the $1^{st}$ 100 time steps.

The number of unknowns was $4\,327,005$. NM for the FEM system is stopped once the relative change in (\ref{eq19}) becomes smaller than $10^{-4}$ and the NM for the net current is stopped once the residuum of the FEM system is reduced by $10^{-5}$.

Experience has shown that a comparison of simulations considering a magnetization curve with measurements fails because of the big influence of hysteresis.


\ifCLASSOPTIONcaptionsoff
  \newpage
\fi



%
\bibliographystyle{IEEEtran}
\bibliography{hollaus}

%








\end{document}